\documentclass[12pt]{amsart}

\usepackage[french, english]{babel}

\newcommand{\baton}[1]{\mathbb #1}
\newcommand{\E}{{\baton E}}
\newcommand{\N}{{\baton N}} 
\newcommand{\R}{{\baton R}}     
\newcommand{\Z}{{\baton Z}}
\newcommand{\CB}{{\mathcal B}}
\newcommand{\CC}{{\mathcal C}} 
\newcommand{\CD}{{\mathcal D}}
\newcommand{\bt}{{\mathbf{t}}}
\newcommand{\bx}{\mathbf{x}}
\newcommand{\bomega}{{\boldsymbol{\omega}}}

\newcommand{\bzero}{{\boldsymbol{0}}}
\newcommand{\one}{{\boldsymbol 1}}
\newcommand{\dis}{\displaystyle}
\newcommand{\norm}[1]{\lVert #1\rVert}
\newcommand{\nnorm}[2]{\Vert #1\Vert_{U^{#2}}}
\DeclareMathOperator{\card}{Card}
\DeclareMathOperator{\ppcm}{\text{\sc ppcm}}
\DeclareMathOperator{\reel}{Re}
\theoremstyle{plain}
\newtheorem{theo}{Th\'eor\`eme}
\newtheorem{prop}{Proposition}
\newtheorem{lemm}{Lemme}

\newtheorem*{ThmSz}{Th\'eor\`eme de Szemer\'edi}

\newtheorem*{ThmSzFini}{Version finie du th\'eor\`eme de Szemer\'{e}di}

\newtheorem*{conj*}{Conjecture}
\theoremstyle{remark}
\newtheorem*{rema*}{Remarque}
\theoremstyle{plain}

\title[Progressions arithm\'etiques et nombres premiers]
{Progressions arithm\'etiques dans les nombres premiers}

\author{ Bernard HOST}

\address{Universit\'e de Marne-La-Vall\'ee \&
UMR-CNRS 8050.
Math\'e\-matiques,
 5 boulevard Descartes,
Champs-sur-Marne.
F--77454 Marne-la-VallŽe C\'edex 2}
\email{bernard.host@univ-mlv.fr}

\begin{document}
\selectlanguage{frenchb}
\maketitle
\centerline{d'apr\`es B. Green et T. Tao}

\begin{abstract}
R\'ecemment, B. Green et T. Tao ont montr\'e que
\emph{l'ensemble des nombres premiers contient des progressions
arithm\'e\-tiques de toutes longueurs},
r\'epondant ainsi \`a une question ancienne \`a la formulation
particuli\`erement simple. La d\'emonstration n'utilise aucune des
m\'ethodes \og transcendantes\fg\ ni aucun des grands th\'eor\`emes de la th\'eorie
analytique des nombres. Elle est \'ecrite
dans un {\it esprit} proche de celui de la th\'eorie ergodique, en
particulier de celui de la preuve par Furstenberg du th\'eor\`eme de
Szemer\'edi, mais  elle n'utilise  aucun th\'eor\`eme provenant de cette
 th\'eorie. La m\'ethode peut ainsi \^{e}tre consid\'er\'ee comme 
\og \'el\'ementaire \fg,
ce qui ne veut pas dire facile.

On se propose de pr\'esenter l'organisation g\'en\'erale de la preuve sans 
 d\'evelopper les calculs.
\end{abstract}
\maketitle

\selectlanguage{english}
\begin{abstract}
B. Green and T. Tao have recently proved that {\em the set of primes
contains arbitrary long arithmetic progressions}, answering to an old
question with a remarkably simple formulation.  The proof does not use
any ``transcendental'' method and any of the deep theorems of analytic
number theory.  It is written in a {\em spirit} close to ergodic
theory and in particular of Furstenberg's proof of Szemer\'edi's
Theorem, but it does not use any result of this theory.  Therefore the
method can be considered as elementary, which does not mean easy.  We
entend here to present the mains ideas of this proof.
\end{abstract} 

\selectlanguage{francais}

\section{Introduction}
\subsection{Le r\'esultat}

Le but de cet expos\'e est de pr\'esenter un travail r\'ecent et
spectaculaire de B. Green et T.~Tao o\`u ils montrent :
\begin{theo}[\cite{GrT}]\label{th:GrT}
L'ensemble des nombres premiers contient des progressions 
arithm\'etiques de toutes longueurs.
\end{theo}
En fait Green et Tao montrent un r\'esultat plus fort : la conclusion 
du th\'eor\`eme reste valable 
si on remplace l'ensemble des nombres premiers par un
sous-ensemble de densit\'e relative positive. De plus, la m\'ethode
employ\'ee permet de d\'eterminer explicitement pour tout $k$ un entier $N$ (tr\`es
grand) tel que l'ensemble des
nombres premiers plus petits que $N$ contienne une progression
arithm\'etique de longueur $k+1$. 

Le th\'eor\`eme~\ref{th:GrT} r\'epond \`a une question fort ancienne bien que difficile \`a
dater exactement. Tr\`es peu de r\'esultats partiels \'etaient connus
jusqu'ici ; citons celui de  van der
Corput~\cite{vdC} qui a montr\'e en 1939 l'existence d'une infinit\'e de
progressions de longueur $3$ dans les nombres premiers.

En 1923, Hardy et Littlewood~\cite{HL}  ont propos\'e 
une conjecture tr\`es g\'en\'erale sur la r\'epartition de certaines
configurations dans les nombres
premiers, qui entra\^\i nerait  une version quantitative pr\'ecise du
th\'eor\`eme~\ref{th:GrT} si elle s'av\'erait exacte.
Ce m\^eme th\'eor\`eme suivrait aussi d'une r\'esolution positive donn\'ee \`a une
conjecture propos\'ee par Erd\"os et Turan~\cite{ET} en 1936 :
\begin{conj*}
Tout sous-ensemble $E$ de $\N^*$ v\'erifiant 
$\dis\sum_{n\in E}\frac 1n=+\infty$ contient des progressions
arithm\'etiques de toutes longueurs.
\end{conj*}
Cette conjecture reste
 totalement ouverte et les m\'ethodes de Green et 
Tao ne permettent pas de s'en approcher. Dans une direction voisine,
Szemer\'edi a montr\'e en 1975 l'existence de progressions sous l'hypoth\`ese plus
forte de la densit\'e positive. Rappelons que 
la densit\'e d'un ensemble d'entiers $E\subset\N$ est :
$$
 d^*(E)=\limsup_{N\to\infty} \frac 1N\card(E\cap[0,N-1])\ .
$$
Le th\'eor\`eme de Szemer\'edi s'\'enonce :
\begin{ThmSz}[\cite{Sz}]
Tout ensemble d'entiers de densit\'e positive contient des progressions 
arithm\'{e}tiques de toutes longueurs.
\end{ThmSz}
Il peut aussi s'exprimer en termes d'ensembles finis d'entiers :
\begin{ThmSzFini}
Pour tout entier $k\geq 2$ et tout r\'eel $\delta>0$ il existe un entier
$N=N(k,\delta)$ tel que tout sous-ensemble $E$ de $[0,N[$ ayant au moins
$\delta N$ \'el\'ements contienne une progression arithm\'etique de longueur
$k+1$.
\end{ThmSzFini}

Ce th\'eor\`eme ne peut \'evidemment pas \^etre utilis\'e directement puisque les 
nombres premiers ont une densit\'e nulle. Cependant il tient une place 
centrale dans la d\'emonstration.

\subsection{La m\'ethode}
Le travail de Green et Tao comporte deux parties tr\`es diff\'erentes.

\medskip
La premi\`ere partie, qui est la plus longue, contient la d\'emonstration 
d'une extension de la version finie du th\'eor\`eme de Szemer\'edi
(th\'eor\`eme~\ref{th:GTSz}).

Dans ce dernier th\'eor\`eme, la quantit\'e $|E|/N\geq\delta$ peut \^etre vue
comme la moyenne sur $[0,N[$ de la fonction indicatrice de $E$. 
L'id\'ee naturelle est de remplacer cette fonction indicatrice par une
fonction nulle en dehors de l'ensemble des nombres premiers, 
mais alors cette fonction ne peut pas \^etre choisie major\'ee par $1$ 
sinon sa moyenne deviendrait arbitrairement petite pour $N$ grand.
Green et Tao montrent un th\'eor\`eme de Szemer\'edi modifi\'e 
(th\'eor\`eme~\ref{th:GTSz}) qui
s'applique \`a une fonction major\'ee par un
\og poids pseudo-al\'eatoire\fg  , c'est \`a dire par une fonction de 
moyenne $1$ dont 
les corr\'{e}lations sont 
voisines de celles qu'on obtiendrait en tirant  au hasard et
ind\'ependamment les valeurs aux points
  $0,1,\dots,N-1$ (section~\ref{subsec:pseudos}). 
Cette utilisation d'une majoration fait penser \`a la
m\'ethode du crible. 

La d\'emonstration de ce \og th\'eor\`eme de Green-Tao Szemer\'edi\fg\
 est \'ecrite dans le
langage des probabilit\'es. Comme tous les espaces de probabilit\'e
sont finis et munis de la mesure uniforme, on pourrait dire qu'elle
utilise seulement  des arguments de d\'enombrement. 
Cette fa\c con de voir serait formellement
correcte mais trop r\'eductrice. En fait la d\'emarche de Green et Tao
s'inspire directement de la th\'eorie ergodique, et plus pr\'{e}cis\'{e}ment de
la d\'emonstration ergodique du th\'eor\`eme de Szemer\'edi donn\'ee par
Furstenberg  (\cite{Fu}, voir aussi~\cite{FKO}). 
Dans les deux cas, le c\oe   ur de la preuve 
est un r\'esultat de d\'ecomposition (proposition~\ref{prop:decomp})
consistant \`a \'ecrire une fonction comme la somme de son esp\'erance
conditionnelle sur une $\sigma$-alg\`ebre bien choisie et d'un reste.
L'esp\'erance conditionnelle est \og liss\'ee\fg\    et dans le cas consid\'er\'e par
Green et Tao elle est m\^eme uniform\'ement born\'ee, ce qui permet d'utiliser le
th\'eor\`eme de Szemer\'edi classique. Le reste se comporte comme une
oscillation al\'eatoire et sa contribution dans les calculs est
n\'egligeable. Les ergodiciens reconna\^{\i}tront la fa\c con dont les
 \og facteurs \fg\  interviennent dans de nombreux probl\`emes. Pour les autres, nous ajoutons que l'article
n'utilise aucun r\'esultat provenant de la th\'eorie ergodique et que sa
lecture ne demande aucune connaissance dans ce domaine.

Cette inspiration ergodique dans une d\'emonstration combinatoire est encore 
plus apparente dans la nouvelle d\'emonstration que T.~Tao vient de donner
du th\'eor\`eme de Szemer\'edi~\cite{Tao}. Nous ne pensons pas que cette
d\'emarche soit artificielle. 
Jusqu'\`a pr\'esent les relations entre ces domaines se r\'esumaient
pratiquement au principe de correspondance de Furstenberg qui permet de 
 montrer, \`a partir de th\'eor\`emes ergodiques, des r\'esultats
combinatoires dont beaucoup n'ont aujourd'hui pas d'autre preuve.
Il appara\^\i t depuis peu des ressemblances de plus en plus prononc\'ees
quoiqu'encore mal comprises entre les objets et 
les m\'ethodes des deux th\'eories. Nous reviendrons dans ces notes sur ce point qui 
m\'erite sans
doute d'\^etre approfondi.

\medskip

Une fois d\'emontr\'e le th\'eor\`eme de Szemer\'{e}di modifi\'e, il reste \`a
construire un poids pseudo-al\'eatoire adapt\'e au probl\`eme pos\'e.
Il s'agit donc ici de th\'eorie des nombres.  
Dans cette partie de l'article~\cite{GrT} les auteurs utilisent une
fonction de von Mangoldt tronqu\'ee et font appel \`a
des outils sophistiqu\'es provenant des travaux de Goldston et  
Y\i ld\i r\i m~\cite{GY} mais,
dans une note non publi\'ee~\cite{Tao2}, T.~Tao explique comment l'argument
peut \^etre modifi\'e pour n'utiliser que les propri\'et\'es les plus \'el\'ementaires des 
nombres premiers et de la fonction $\zeta$.
 C'est cette approche que nous suivons ici
 en nous inspirant de notes manuscrites de J.-C.~Yoccoz.

\medskip
Dans cet  expos\'e, qui ne contient aucune d\'emonstration compl\`ete,
on se propose de pr\'esenter de fa\c con assez d\'etaill\'ee
l'organisation  de la preuve et de  donner une id\'ee des
m\'ethodes employ\'ees \`a chaque \'etape.
Le lecteur press\'e pourra se  limiter \`a la section~\ref{sec:GrTSz} qui 
contient la formulation pr\'ecise des d\'efinitions et r\'esultats
correspondant aux deux grandes parties auxquelles on vient de faire
allusion, encore que
la d\'efinition des normes de Gowers
(sous-sections~\ref{subsec:def-gowers} et~\ref{subsec:comment})
ait son int\'er\^et propre. 
Le r\'esultat  de d\'ecomposition (proposition~\ref{prop:decomp}) est \'enonc\'e dans la
sous-section~\ref{subsec:sigma} et montr\'e dans la
section~\ref{sec:decomp}.
La deuxi\`eme partie de la preuve, c'est \`a dire la construction du
poids pseudo-al\'eatoire, est contenue dans la
section~\ref{sec:pseudo}.

\subsection{Conventions et notations}

Quand $f$ est une fonction d\'efinie sur un ensemble fini $A$,
l'\emph{esp\'erance} de $f$ sur $A$, not\'ee 
$\E(f(x)\mid x\in A)$ ou $\E(f\mid A)$, est la moyenne arithm\'etique de $f$ sur 
$A$ ; la m\^eme est utilis\'ee pour les fonctions de 
plusieurs variables.

Dans  toute la suite, $k\geq 2$ est un entier que nous consid\'erons
comme une constante. L'objectif est de montrer l'existence d'une
progression arithm\'{e}tique de longueur $k+1$ dans les nombres premiers.
La progression est cherch\'ee dans l'intervalle 
$[0,N[$,  o\`u  $N$ est un (grand) entier 
qu'il est souvent n\'ecessaire de supposer premier. 
On identifie $[0,N[$ au groupe $\Z_N=\Z/N\Z$. 

Il est crucial dans la preuve de contr\^oler la mani\`ere dont toutes les 
estimations d\'ependent de $N$ et  nous adoptons les conventions suivantes.  
 Dans chaque \'enonc\'e $N$ est suppos\'e fix\'e mais toutes les constantes sont 
ind\'ependantes de $N$.
Nous notons $o(1)$ une quantit\'e tendant vers $0$ quand $N$ tend vers
l'infini, uniform\'ement par rapport \`a tous les param\`etres sauf
\'eventuellement ceux not\'es en indice. La notation $O(1)$ est employ\'ee
avec une signification similaire.

\section{Poids
pseudo-al\'eatoires\\ et th\'eor\`eme de Green-Tao Szemer\'edi\\ }
\label{sec:GrTSz}

Green et Tao g\'en\'eralisent une formulation classique du th\'eor\`eme de
Szemer\'edi,
qui est celle sous laquelle Gowers~\cite{G}  l'a red\'emontr\'e r\'ecemment.

\begin{theo}
\label{prop:SzInt}
Pour tout  r\'eel $\delta>0$ il existe 
un constante $c(\delta)>0$ tel que, pour toute fonction $f\colon\Z_N\to\R$ avec 
$$
0\leq f(x)\leq 1\text{ pour tout $x$  et }\E\bigl(f\mid \Z_N\bigr)
\geq\delta
$$
on ait
\begin{equation}
\label{eq:SzInt}
\E\bigl(f(x)f(x+t)\dots f(x+kt)\mid x,t\in\Z_N\bigr)\geq 
c(\delta)\ .
\end{equation}
\end{theo} 
La version finie du th\'eor\`eme de Szemer\'edi se d\'eduit de ce th\'eor\`eme en 
prenant pour $f$ la fonction indicatrice d'un sous-ensemble de $[0,N[$. 
 Green
et Tao s'affranchissent de la condition $f\leq 1$ en la rempla\c cant
par l'hypoth\`ese que $f$ est major\'ee par 
un \emph{poids pseudo-al\'eatoire} ;
cette notion  sera d\'efinie plus loin.

\subsection{Les deux ingr\'edients de la preuve du 
th\'eor\`eme~\ref{th:GrT}}

 Nous appelons \og th\'eor\`eme de Green-Tao Szemer\'edi\fg\   l'extension
suivante du th\'eor\`eme de Szemer\'{e}di : 
\begin{theo}
\label{th:GTSz} 
Soit $\nu\colon\Z_N\to\R^+$ un poids pseudo-al\'eatoire 
(voir la sous-section~\ref{subsec:pseudos}). 
 Pour tout r\'eel $\delta>0$ il existe 
une constante $c'(\delta)>0$ satisfaisant la propri\'et\'e suivante. 
Pour toute fonction $f\colon\Z_N\to\R$  telle que
$$
0\leq f(x)\leq\nu(x)\text{ pour tout $x$ et } 
\E\bigl(f\mid \Z_N\bigr)\geq\delta
$$
on a
\begin{equation}
\label{eq:GTSz}
\E\bigl(f(x)f(x+t)\dots f(x+kt)\mid x,t\in\Z_N\bigr)\geq 
c'(\delta)-o(1)\ .
\end{equation}
\end{theo}
La d\'emonstration de ce th\'eor\`eme, qui occupe une part importante de
l'article de Green et Tao, est r\'esum\'ee dans les
sections~\ref{sec:normes} et~\ref{sec:decomp}. Pour
l'appliquer aux nombres premiers, il faut une fonction $f$ et un poids
$\nu$
convenables dont l'existence est donn\'ee par le th\'eor\`eme suivant. 
\begin{theo}
\label{th:existe-nu}
Il existe une constante positive $\delta$, un  poids pseudo-al\'eatoire 
$\nu\colon\Z_N\to\R^+$ et une
fonction $f\colon\Z_N\to\R$ avec
\begin{gather*}
f \text{ est nulle en dehors de l'ensemble des nombres premiers ;}\\
0\leq f(x)\leq\nu(x)\text{ pour tout $x$ ;}\\ 
\E(f\mid\Z_N)\geq\delta\  ; \\
\norm f_{L^\infty}= O(1)\,\log N\ .
\end{gather*}
\end{theo}

La construction de $f$ et $\nu$ est faite dans la
section~\ref{sec:pseudo}. Nous montrons maintenant comment le th\'eor\`eme~\ref{th:GrT} 
d\'ecoule des th\'eor\`emes~\ref{th:GTSz} et~\ref{th:existe-nu}. 

\begin{proof}
Soient $\delta$, $f$ et $\nu$ comme dans le th\'eor\`eme~\ref{th:existe-nu}. 
Il existe un intervalle $J\subset[0,N[$, de longueur plus petite que
$ N/2$ et
tel que $\E(\one_Jf\mid\Z_N)\geq \delta/3$. Nous utilisons le
th\'eor\`eme~\ref{th:GTSz} avec la fonction $f$ remplac\'ee par $\one_Jf$ 
et le r\'eel $\delta$
remplac\'e par $\delta/3$. 

La contribution dans l'esp\'erance~\eqref{eq:GTSz} 
des termes o\`u $t=0$ est major\'ee par 
$N^{-1}\norm f_{L^\infty}^{k+1}= o(1)$, et est donc inf\'erieure \`a
$c'(\delta)$ si $N$ est assez grand.  Il existe donc dans ce cas
 $x,t\in\Z_N$ avec $t\neq 0$ tels que 
$f(x)f(x+t)\dots f(x+kt)\neq 0$.
Rappelons que dans cette
expression $x, x+t,\dots,x+kt$ sont consid\'er\'es comme des \'el\'ements de
$\Z_N$ et que donc l'addition est modulo $N$. Si nous consid\'{e}rons $x$
et $t$ comme des entiers appartenant \`a l'intervalle $[0,N[$ 
nous obtenons que $f$ est non nulle aux points  
$x,x+t\bmod N,\dots,x+kt\bmod N$. 
Comme elle est nulle en dehors de l'intervalle $J$ de longueur 
$<N/2$,
tous ces
entiers  appartiennent \`a cet intervalle et on en d\'eduit facilement
qu'ils forment une progression arithm\'etique non triviale de longueur
$k+1$. Enfin, $f$ est nulle en dehors de l'ensemble des nombres
premiers et on a bien une progression form\'ee de nombres premiers.
\end{proof}

\subsection{D\'efinition des poids pseudo-al\'eatoires}
\label{subsec:pseudos}
Dans les th\'eor\`emes pr\'ec\'edents nous avons consid\'er\'e un
\emph{poids pseudo-al\'eatoire} comme une fonction d\'efinie sur
$\Z_N$. 
Il s'agit plus pr\'{e}cis\'{e}ment de la donn\'ee, pour chaque nombre premier $N$, d'une
fonction $\nu=\nu_N\colon \Z_N\to\R^+$, de sorte que soient
satisfaites deux conditions
asymptotiques appel\'ees 
 \emph{condition sur les formes lin\'eaires} et  \emph{condition sur les 
corr\'elations}.

\subsubsection*{La condition sur les formes lin\'eaires.} 
Ici $m_0,t$ et $L$ sont des constantes enti\`eres (ne d\'ependant que de $k$) 
que nous n'explicitons pas.\\
Soient $m\leq m_0$ un entier et $\psi_1,\dots,\psi_m$   des 
applications de $\Z_N^t$ dans $\Z_N$ de la forme 
\begin{equation}
\label{eq:def-psi-i}
\psi_i(\bx)=b_i+\,\sum_{j=1}^t L_{i,j}x_j
\end{equation}
o\`u $\bx=(x_1,\dots,x_t)$ et 
\begin{itemize}
\item
pour tout $i$, $b_i$ est un entier; 
\item
pour tous $i,j$, $L_{i,j}$ est un entier avec $|L_{i,j}|\leq L$;
\item
aucun des vecteurs $(L_{i,j})_{1\leq j\leq t}\in\Z^t$ 
n'est nul
et ces vecteurs sont deux \`a deux non colin\'eaires
\end{itemize}
alors la condition sur les formes lin\'eaires stipule que
\begin{equation}
\label{eq:lineaire}
\E\Bigl(\nu(\psi_1(\bx))\dots\nu(\psi_m(\bx))\,\big|\,
\bx\in\Z_N^t\Bigr)=1+o(1)\ .
\end{equation}
 D'apr\`es nos conventions, 
la quantit\'e $o(1)$ tend vers $0$ quand $N$ tend vers l'infini
ind\'{e}pendamment du choix des fonctions $\psi_i$ et en 
particulier du choix des $b_i$, qui ne sont pas suppos\'es born\'es.
Remarquons que la condition des formes lin\'eaires  entra\^\i ne que la m\^eme  majoration reste 
valable s'il y a moins de $t$ variables.  En particulier,
$$
 \E(\nu\mid\Z_N)=1+o(1)\ .
$$

\subsubsection*{La condition des corr\'elations.}
Ici encore, $q_0$ est une constante enti\`ere  que nous 
n'explicitons pas.\\
La condition sur les corr\'elations stipule 
qu'il existe une fonction $\tau\colon\Z_N\to\R^+$ avec
$$
\text{pour tout $p\geq 1$, }\E(\tau^p(x)\mid x\in\Z_N)=O_p(1) 
$$ 
telle que,  pour tout $q\leq q_0$ et tous $h_1,\dots,h_q\in\Z_N$,
 distincts ou confondus, on ait
\begin{equation}
\label{eq:correlation}
\E\bigl( \nu(x+h_1) \nu(x+h_2) \dots \nu(x+h_q)\mid 
x\in\Z_N\bigr)
\leq\sum_{1\leq i\leq j\leq q}\tau(h_i-h_j)\ .
\end{equation}
Nous remarquons que, si $\nu$ est un poids pseudo-al\'eatoire, alors $(1+\nu)/2$ 
en est \'egalement un.  On  
peut donc sans perte de g\'en\'eralit\'e se restreindre au 
cas o\`u $\nu(x)>0$ pour tout $x$.

\section{Les normes de Gowers}
\label{sec:normes}
\subsection{La d\'efinition}
\label{subsec:def-gowers}
Il y a quelques ann\'ees Gowers a propos\'e une nouvelle preuve 
du th\'eor\`eme de Szemer\'edi~\cite{G} \`a base d'analyse harmonique et de
combinatoire. Dans sa d\'emonstration il 
 a introduit une suite de normes sur l'espace  $\CC(\Z_N)$ des 
fonctions sur $\Z_N$ \`a valeurs r\'eelles et les a utilis\'{e}es pour 
contr\^oler les esp\'erances qui apparaissent dans le th\'eor\`eme~\ref{prop:SzInt}. 
Green et Tao les utilisent \'egalement pour contr\^oler les esp\'erances du 
th\'eor\`eme~\ref{th:GTSz}. Nous donnons ici leur d\'efinition.

Pour $f\in\CC(\Z_N)$ on  d\'efinit par r\'ecurrence
les quantit\'es
$\nnorm f d$\,, $d\geq 1$, par 
\begin{align*}
\nnorm f 1
&=\bigl|\E(f\mid \Z_N)\Bigr|\ ;\\
\nnorm f{d+1}
&=\Bigl(\E\bigl(\nnorm{f\cdot f_t}d^{2^d} \mid
t\in\Z_N\bigr)\Bigr)^{1/2^{d+1}}
\text{ pour }d\geq 1
\end{align*}
o\`u $f_t$ est la fonction  $x\mapsto f(x+t)$.

Ces quantit\'es peuvent \^etre aussi donn\'ees par une formule close. 

Pour $\bomega=(\omega_1,\dots,\omega_d)\in\{0,1\}^d$ et 
$\bt=(t_1,\dots,t_d)\in\Z_N^d$ notons 
$$
\bomega\cdot\bt=\omega_1 
t_1+\omega_2 t_2+\dots+\omega_d t_d\ .
$$
 On a alors
\begin{equation}
\label{eq:def-normes}
\nnorm fd
=\Bigl(\E\Bigl(\prod_{\bomega\in\{0,1\}^d}f(x+\bomega\cdot\bt\,\big|\,
x\in\Z_N,\; \bt\in\Z_N^d\Bigr)\Bigr)^{1/2^d}\ .
\end{equation}
On v\'erifie alors facilement que $\nnorm f2$ est la norme $\ell^4$ de 
la transform\'{e}e de Fourier de $f$ et  que $\nnorm f{d+1}\geq \nnorm 
fd$ pour tout $d\geq 1$. De plus on obtient :

\begin{prop}[In\'egalit\'e de Cauchy-Schwarz-Gowers]
\label{prop:CSG}\strut
Si $f_\bomega$, $\bomega\in\{0,1\}^d$, sont $2^d$ fonctions r\'eelles sur 
$\Z_N$ on a
\begin{equation}
\label{eq:CSG}
\Bigl|\E\Bigl(\prod_{\bomega\in\{0,1\}^d}f_\bomega(x+\bomega\cdot\bt\,\big|\,
x\in\Z_N,\; \bt\in\Z_N^d\Bigr)\Bigr|\leq
\prod_{\bomega\in\{0,1\}^d}\nnorm {f_\omega}d\ .
\end{equation}
\end{prop}

On en d\'eduit :
\begin{prop}
Pour $d\geq 2$ l'application $f\mapsto\nnorm fd$ est 
une norme sur $\CC(\Z_N)$.
\end{prop}
On peut facilement \'etendre ces d\'efinitions au cas des fonctions \`a
valeurs complexes.
\subsection{Commentaires}
\label{subsec:comment}
 Pour $d>2$ la norme  $\nnorm\cdot d$ est
assez difficile \`a 
interpr\'eter car elle ne peut apparemment pas \^etre exprim\'ee 
 au moyen des normes classiques. La d\'efinition n'est pas simplifi\'ee 
par l'usage de la transform\'ee de Fourier ; par exemple, la norme 
$\nnorm\cdot 3$ d'une fonction est la m\^eme que celle de sa 
transform\'ee de Fourier (\`a une normalisation pr\`es).  

Cette difficult\'e provient sans doute du fait que ces normes ont 
un aspect non commutatif.
En effet, il est clairement possible de d\'efinir des normes similaires sur
l'espace $\CC_K(G)$ des fonctions continues \`a support compact sur un groupe
ab\'elien localement compact $G$. Mais il est sans doute moins \'evident 
que la d\'efinition de 
 la norme $\nnorm\cdot d$ s'\'etend au cas o\`u
$G$  est localement compact nilpotent d'ordre $d-1$, et qu'elle peut
m\^eme \^etre d\'efinie sur
 $\CC_K(G/\Gamma)$ lorsque $\Gamma$ est un sous-groupe ferm\'e 
d'un groupe $G$ de ce type.

D'une mani\`ere ind\'ependante,  des semi-normes 
$\lvert\!|\!| \cdot|\!|\!\rvert_d$\,, $d\geq 
1$, ont r\'ecemment \'et\'e introduites~\cite{HK} en 
th\'eorie ergodique dans l'\'etude de questions relatives au th\'eor\`eme de 
Szemer\'edi o\`u elles servent \'egalement \`a contr\^oler des esp\'erances
ressemblant \`a celles du th\'eor\`eme~\ref{prop:SzInt}. 
La d\'efinition de ces semi-normes, nettement plus 
compliqu\'ee, ne sera pas donn\'{e}e ici mais elle est formellement assez 
similaire \`a celle des normes de Gowers. Ces semi-normes ont une 
interpr\'etation simple : elles sont li\'ees \`a l'existence de 
quotients du syst\`eme munis d'une structure   
d'espace homog\`ene d'un groupe de Lie nilpotent. 

Si on admet que les 
ressemblances de plus en plus nombreuses qui apparaissent entre les deux 
th\'eories ne sont pas fortuites, il est alors possible de conjecturer 
que les normes de Gowers s'interpr\`etent au moyen d'une sorte de 
transform\'ee de Fourier nilpotente, m\^eme lorsque le groupe est ab\'elien.

\subsection{Normes de Gowers et progressions arithm\'etiques}
La proposition suivante g\'en\'eralise un r\'esultat analogue de Gowers
\'etabli sous l'hypoth\`ese plus forte que toutes les fonctions sont
born\'ees par $1$.  Sa d\'emonstration consiste en une suite ing\'enieuse
d'applications de l'in\'egalit\'e de Cauchy-Schwarz, de changements de
variables et de la condition sur les formes lin\'eaires.
\begin{prop}
\label{prop:normes-progs}
Soient $\nu$ un poids pseudo-al\'eatoire et $f_0,f_1,\dots,f_{k}$ 
des fonctions sur $\Z_N$ v\'erifiant 
$$
 |f_j(x)|\leq 1+\nu(x)\text{ pour tout }x\in\Z_N\text{ et tout $j$ avec }
0\leq j\leq k\ .
$$
Alors
$$
\Bigl| \E\Bigl(\prod_{j=0}^{k}f_j(x+jt)\,\big|\, x,t\in\Z_N\Bigr)
\Bigr|\leq 2^{k+1}\min_{0\leq j\leq k}\nnorm{f_j}{k} +o(1)\ .
$$
\end{prop}
L'utilisation que font Green et Tao de cette proposition est tr\`es
diff\'erente de la mani\`ere dont Gowers utilise le r\'{e}sultat analogue pour 
les fonctions born\'ees.

Ce dernier proc\`ede par dichotomie.\\
 Soient $f$
une fonction sur $\Z_N$ et $c=\E(f\mid\Z_N)$. 
Si $\nnorm{f-c}{k}$
est petit, alors l'esp\'erance~\eqref{eq:SzInt} est peu diff\'erente de
l'esp\'erance obtenue en rempla\c cant $f$ par $c$ et elle est donc grande. Si au
contraire cette norme est grande, alors Gowers montre que la
restriction de $f$ \`a un sous-ensemble pas trop petit de $\Z_N$
pr\'esente des r\'egularit\'es qui sont ensuite exploit\'ees.

Green et Tao utilisent une d\'ecomposition, o\`u les normes de Gowers
jouent un r\^ole tr\`es proche de celui jou\'e par les
 semi-normes $\lvert\!|\!| \cdot|\!|\!\rvert_d$ dans~\cite{HK}.
Lorsque  $f$ est une fonction major\'ee par un
poids pseudo-al\'eatoire, elle  peut s'\'ecrire (essentiellement)
comme la somme
d'une fonction ayant une petite norme et d'une fonction born\'ee
 qui est son esp\'erance conditionnelle
sur une $\sigma$-alg\`ebre (proposition~\ref{prop:decomp}). La proposition~\ref{prop:normes-progs}
permet alors de  borner la contribution provenant de la fonction de petite norme.
En th\'eorie ergodique on \'ecrit chaque fonction comme somme  de son
esp\'erance conditionnelle sur une $\sigma$-alg\`ebre adapt\'ee  et
d'une fonction de semi-norme nulle. On utilise ensuite le fait que
cette $\sigma$-alg\`ebre a une interpr\'etation \og 
g\'eom\'etrique\fg\   assez simple.

\subsection{$\sigma$-alg\`ebres sur $\Z_N$ et un r\'esultat de 
d\'ecomposition} 
\label{subsec:sigma}
$\Z_N$ \'etant fini, toute  $\sigma$-alg\`ebre $\CB$ sur $\Z_N$ 
est  d\'efinie par une partition de cet
ensemble : les \'el\'ements de $\CB$ sont les r\'eunions d'atomes de cette
partition et les fonctions $\CB$-mesurables sont les fonctions
constantes sur chaque atome. Quand $f$ est une fonction sur $\Z_N$, son
esp\'erance conditionnelle par rapport \`a $\CB$ est la fonction
$\CB$-mesurable d\'efinie par
\begin{multline*}
\text{si  $A$ est l'atome de $\CB$ contenant $x$,}\\
 \E(f\mid\CB)(x)=\E\bigl(f\mid  A\bigr)
=\frac{\E(\one_A\,f\mid\Z_N)}{\E(\one_A\mid\Z_N)}\ .
\end{multline*}

La proposition suivante est la cl\'e de la d\'emonstration du
th\'eor\`eme~\ref{th:GTSz}.	

\begin{prop}
\label{prop:decomp} 
Soit $\nu$ un poids pseudo-al\'eatoire.
Pour tout $\varepsilon>0$ il existe un entier $N_0(\varepsilon)$ tel que
pour tout $N>N_0(\varepsilon)$ on ait la propri\'et\'e suivante.

Soit $f$ une fonction sur $\Z_N$ avec $0\leq f(x)\leq\nu(x)$ 
pour tout $x$.
 Alors il existe une $\sigma$-alg\`ebre $\CB$ sur $\Z_N$, un
sous-ensemble $\Omega$ de $\Z_N$ appartenant \`a $\CB$ avec
\begin{gather}
\label{eq:nu-sur-omega}
 \E(\nu\cdot\one_\Omega\mid\Z_N)=o_\varepsilon(1)\ ;\\
\label{eq:norme-infty}
 \bigl\Vert(1-\one_\Omega)\cdot\E(\nu-1\mid\CB)\bigr\Vert_{L^\infty}=o_\varepsilon(1)\
\ ;\\
\label{eq:f-esp-f}
\bigl\Vert(1-\one_\Omega)\cdot(f-\E(f\mid\CB))\bigr\Vert_{U^k}\leq\varepsilon
\ .
\end{gather}
\end{prop}

On donne dans la  section~\ref{sec:decomp} un r\'esum\'e de la preuve de
cette proposition.
Admettant ce r\'esultat pour le moment, nous indiquons comment on peut en
d\'eduire le th\'eor\`eme~\ref{th:GTSz}. 

\subsection{D\'emonstration du th\'eor\`eme~\ref{th:GTSz} \`a partir des
propositions~\ref{prop:CSG} et~\ref{prop:decomp}}
Soient $\nu, f$ et $\delta$ comme dans le th\'eor\`eme. Soient $\varepsilon>0$
un param\`etre suffisamment petit et  $\CB$, $\Omega$ comme dans la
proposition~\ref{prop:decomp}. Nous supposons que $N$ est suffisamment
grand. Posons
$$
g=(1-\one_\Omega)\cdot\E(f\mid\CB) 
\text{ et }
h=(1-\one_\Omega)\cdot\bigl(f-\E(f\mid\CB)\bigr)\ .
$$
Comme $f\leq\nu$ nous avons
\begin{multline}\label{eq:minor-g}
\E(g\mid\Z_N)
 \geq\E(f\mid\Z_N)-\E\bigl(\one_\Omega\cdot 
\E(\nu\mid\CB)\mid\Z_N\bigr)\\
=\E(f\mid\Z_N)-\E(\one_\Omega\cdot\nu\mid\Z_N) \geq
\delta- o_\varepsilon(1)
\end{multline}
car $\Omega\in\CB$ et d'apr\`es~\eqref{eq:nu-sur-omega}. De plus
\begin{equation}
\label{eq:borne-g}
 0\leq g\leq (1-\one_\Omega)\cdot\E(\nu\mid\CB) \leq 1+o_\varepsilon(1)
\end{equation}
d'apr\`es~\eqref{eq:norme-infty}. Ainsi, $|h|\leq f+g\leq
1+\nu+o_\varepsilon(1)$ et par ailleurs $\nnorm h k\leq \varepsilon$
d'apr\`es~\eqref{eq:f-esp-f}.

Comme $0\leq g+h\leq f$, l'esp\'erance~\eqref{eq:GTSz}
 apparaissant dans le th\'eor\`eme est minor\'ee par la m\^eme esp\'erance avec
$f$ remplac\'ee par $g+h$. Cette derni\`ere expression s'\'ecrit comme somme
de $2^{k+1}$ esp\'erances de la forme
\begin{equation}
\label{eq:f-i}
\E\bigl(f_0(x)f_1(x+t)\dots f_k(x+kt)\mid x,t\in\Z_N\Bigr)
\end{equation}
o\`u chacune des fonctions $f_i$, $0\leq i\leq k$, est \'egale \`a $g$ ou \`a
$h$. Ainsi, $|f_i|\leq 1+\nu+o_\varepsilon(1)$ pour tout $i$.

Le terme principal est celui o\`u toutes les fonctions $f_i$ sont \'egales
\`a $g$; en effet la majoration~\eqref{eq:borne-g} permet d'utiliser
le th\'eor\`eme de Szemer\'edi (th\'eor\`eme~\ref{prop:SzInt}) et la
minoration~\eqref{eq:minor-g} entra\^ine donc que ce terme  est  minor\'e 
par $c\bigl(\delta-o_\varepsilon(1)\bigr)$. Tous les autres
termes ont une valeur absolue major\'ee par $2^{k+1}\varepsilon+o_\varepsilon(1)$ d'apr\`es la 
proposition~\ref{prop:normes-progs}. 

En choisissant $\varepsilon$ assez petit nous obtenons donc la
minoration annonc\'ee de l'esp\'erance~\eqref{eq:GTSz}, avec
$c'(\delta)=c(\delta)$.\hfill\qed

\begin{rema*}
Le r\'esultat obtenu est plus fort que ce qui est r\'eellement n\'ecessaire,
\`a savoir $c'(\delta)>0$. Il serait sans doute possible de modifier la 
d\'emonstration en affaiblissant les
conditions impos\'ees \`a $\nu$ tout en conservant la propri\'et\'e annonc\'ee.
\end{rema*}

\section{D\'emonstration de la proposition 3.4}
\label{sec:decomp}
Cette section est la plus technique de ces notes et les lecteurs qui 
ne seraient pas
int\'eress\'es par les d\'etails sont invit\'es \`a passer directement \`a la
suivante.

\subsection{Les fonctions duales}

Soit $f$ une fonction r\'eelle sur $\Z_N$.
 Pour $x\in\Z_N$ d\'efinissons
$$
 \CD f(x)=\E\Bigl(\prod_{\substack{\bomega\in\{0,1\}^{k}\\
\bomega\neq\bzero}} f(x+\bomega\cdot\bt)\,\Big|\,\bt\in\Z_N^{k}\Bigr)
$$
o\`u $\bzero$ repr\'esente l'\'el\'ement $(0,0,\dots,0)$ de $\{0,1\}^{k}$.
$\CD f$ est appel\'ee la \emph{fonction duale} (d'ordre $k$) de $f$. 

\'Ecrivons $\langle\cdot,\cdot\rangle$ le produit scalaire sur
$\CC(\Z_N)$ donn\'e par $\langle f,g\rangle=\E(fg\mid\Z_N)$. La
d\'efinition des normes et l'in\'egalit\'e de Cauchy-Schwarz-Gowers
entra\^{\i}nent imm\'ediatement :
\begin{lemm}
\label{lem:dualite}
Pour toute fonction $f \in\CC(\Z_N)$,
$$
 \nnorm f {k}=\langle f,\CD f\rangle =
\sup \Bigl\{ \bigl|\langle f,\CD g\rangle\bigr|\;;\;
g\in\CC(\Z_N),\; \nnorm g {k}\leq 1\Bigr\}\ . 
$$
\end{lemm}
Ainsi, la boule unit\'e pour la norme duale  de $\nnorm\cdot k$ est
l'enveloppe convexe de l'ensemble 
$\{ \CD f\;;\; \nnorm fk\leq 1\}$.
Cette norme duale n'est malheureusement pas une norme d'alg\`ebre (la
norme d'un produit n'est pas major\'ee par le produit des normes), ce
qui simplifierait beaucoup la d\'emonstration. Dans sa preuve du
th\'eor\`eme de Szemer\'edi~\cite{Tao}, Tao construit une norme d'alg\`ebre
qui est major\'ee par la norme duale. Cette construction est
formellement tr\`es proche de la construction de la \og tour d'extensions 
isom\'etriques\fg\ de Furstenberg. 

Pour comprendre le r\^ole jou\'e par les fonctions duales, imaginons la
situation o\`u nous avons une fonction $f$ telle que  $\nnorm f k$ 
soit \og grande \fg\   et que $\norm{\CD f}_{L^2}$ ne soit pas 
\og trop grande\fg.
Supposons aussi que nous savons construire une  $\sigma$-alg\`ebre
$\CB$ par rapport \`a laquelle $\CD f$ est mesurable au 
moins approximativement. Alors, comme le produit scalaire de $f$ et
$\CD f$ est grand, l'esp\'erance  $\E(f\mid\CB)$ aura une norme $L^2$
assez grande. Cette m\'ethode est utilis\'ee de mani\`ere it\'erative 
dans les
sous-sections suivantes pour construire la $\sigma$-alg\`ebre de la 
proposition~\ref{prop:decomp}.

\subsection{Poids pseudo-al\'eatoires et fonctions duales}

Dans toute la suite de cette section, $\nu$ d\'esigne un poids
pseudo-al\'eatoire et nous \'etudions les propri\'et\'es des fonctions duales 
des fonctions major\'ees par $\nu$ ou par $1+\nu$. 

Rappelons que la condition sur les formes lin\'eaires entra\^\i ne que
 $\E(\nu\mid \Z_N)=1+o(1)$. On a  plus pr\'{e}cis\'{e}ment
\begin{equation}
\label{eq:norm-nu}
\nnorm{ \nu -1}k=o(1)\ .
\end{equation}
Nous obtenons de m\^eme :
\begin{lemm}
\label{lem:Df-bornee}
Si $f$ est une fonction sur $\Z_N$ v\'erifiant $|f|\leq 1+\nu$ 
alors $\norm{\CD f}_{L^\infty}\leq 2^{{2^k}-1}+o(1)$.
\end{lemm}
Nous notons d\'esormais $I$ un intervalle ferm\'e born\'e de $\R$ tel que 
$\CD f(x)\in I$ pour tout $x$ et toute fonction $f$ avec $|f|\leq
1+\nu$.

\begin{prop}
\label{prop:fonctions}
Soient $m\geq 1$ un entier et $f_1,\dots,f_m$ des fonctions sur $\Z_N$
v\'erifiant $|f_i|\leq 1+\nu$ pour tout  $i$, et soit
$\Phi$ une fonction continue sur le cube $I^m$.  Alors la 
fonction $\psi$ sur $\Z_N$ d\'efinie par
$$
 \psi(x)=\Phi\bigl(\CD f_1(x),\dots,\CD f_m(x)\bigr)
$$
satisfait la relation 
$$
 \langle \nu-1, \psi\rangle=o_{m,\Phi}(1)\ .
$$
De plus, cette estimation est uniforme en $\Phi$ si l'on impose \`a
cette fonction de rester dans un compact au sens de la convergence
uniforme.
\end{prop}
Pour montrer cette proposition on se ram\`ene facilement au cas o\`u 
$\Phi(x_1,x_2,\dots,x_m)=x_1x_2\ldots x_m$ 
et on utilise la condition des corr\'elations. 
C'est le seul endroit de la preuve o\`u cette
condition est utilis\'ee.  

\subsection{Construction d'une $\sigma$-alg\`ebre}
Nous introduisons ici une cons\-truction qui sera utilis\'ee de mani\`ere
r\'ep\'et\'ee dans la section suivante pour montrer la proposition~\ref{prop:decomp}. Ici
$\varepsilon>0$ est un param\`etre et $\sigma\in]0,1/2[$ est un param\`etre 
accessoire qui devra \^etre choisi soigneusement en fonction de
$\varepsilon$. On se donne une fonction continue $\psi\colon\R\to[0,1]$, \`a support
dans $[0,1]$ et \'egale \`a $1$ sur $[\sigma,1-\sigma]$. 
Nous supposons toujours que $N$ est suffisamment grand.

Soit $f$ une fonction sur $\Z_N$ avec $|f|\leq 1+\nu$ et notons 
$F=\CD f$.

Soient $\alpha\in]0,1]$ et $\CB$ la $\sigma$-alg\`ebre dont les atomes
sont les ensembles $A$ de la forme
\begin{equation}
\label{eq:def-atome}
 A=\{x\in\Z_N\;;\; \varepsilon^{2^{k+1}}(n+\alpha)\leq
F(x)<\varepsilon^{2^{k+1}}(n+1+\alpha)\}
\end{equation}
o\`u $n$ est un entier tel que cet ensemble ne soit pas vide.
Le param\`etre $\alpha$ est introduit pour \'eviter les effets de bord : il
pourrait en effet arriver que les valeurs de la fonction $F$ s'accumulent
pr\`es des points $n\varepsilon^{2^{k+1}}$ mais, pour un choix convenable de 
$\alpha$, l'ensemble 
$$
 E=\bigcup_{n\in\Z}\bigl\{x\in\Z_N\;;\; \varepsilon^{2^{k+1}}(n+\alpha-\sigma)\leq
 F(x)\leq
 \varepsilon^{2^{k+1}}(n+\alpha+\sigma)\bigr\}
$$
v\'erifie 
\begin{equation}
\label{eq:choix-alpha}
 \E\bigl(\one_E\cdot(1+\nu)\mid\Z_N\bigr)=\sigma\,O(1)\ .
\end{equation}
Par construction, on a clairement,
\begin{equation}
\label{eq:appro-infty}
\bigl\Vert F-\E(F\mid\CB)\bigr\Vert_{L^\infty}\leq\varepsilon^{2^{k+1}}
\end{equation}
Comme $F$ est born\'ee (lemme~\ref{lem:Df-bornee}), le nombre d'atomes de $\CB$ est un
$O_{\varepsilon}(1)$.

Appelons un atome $A$ de $\CB$ \emph{mauvais} si
$\E\bigl((1+\nu)\,\one_A\mid\Z_N\bigr)<\sigma^{1/2}$ et notons $\Omega$ 
la r\'eunion des mauvais atomes. Alors $\Omega\in\CB$ et 
\begin{equation}
\label{eq:int-omega}
\E\bigl((1+\nu)\,\one_\Omega\mid\Z_N\bigr)=\sigma^{1/2}\,O_\varepsilon(1)
\ .
\end{equation}

Soient maintenant $A$ un  \emph{bon} atome, $n$ l'entier
correspondant dans la d\'efinition~\eqref{eq:def-atome} et 
$J= [\varepsilon^{2^{k+1}}(n+\alpha,\varepsilon^{2^{k+1}}(n+1+\alpha)[$ .
Posons  $\Phi_A(x)=\psi\bigl(\varepsilon^{-2^{k+1}}(x-n-\alpha)\bigr)$. 
La proposition~\ref{prop:fonctions} permet de majorer $\Phi_A\circ F$
et la propri\'et\'e~\eqref{eq:choix-alpha} permet de contr\^oler le terme
d'erreur $\one_A-\Phi_A\circ F=(\one_J-\Phi_A)\circ F$. Nous obtenons
$$
 \E\bigl((\nu(x)-1)\one_A(x)\mid x\in A\bigr)
=\sigma^{1/2}\, O_\varepsilon(1)+o_{\varepsilon,\sigma}(1)\ .
$$
On en d\'eduit que 
\begin{equation}
\label{eq:borne-unif}
\bigl\Vert(1-\one_\Omega)\cdot \E(\nu-1\mid\CB) \bigr\Vert_{L^\infty}
=0_{\varepsilon}(1)\ .
\end{equation}

Supposons maintenant qu'au lieu d'une fonction $f$ nous avons une
famille finie $(f_1,\dots,f_m)$ de fonctions v\'erifiant toutes
$|f_i|\leq 1+\nu$. Alors par la m\^eme m\'ethode nous pouvons construire
une $\sigma$-alg\`ebre $\CB$ et un ensemble $\Omega\in\CB$ v\'erifiant 
les propri\'et\'es~\eqref{eq:int-omega} et~\eqref{eq:borne-unif} et tels
que l'approximation uniforme~\eqref{eq:appro-infty} soit valable pour
chacune des fonctions $F_i=\CD f_i$, tous les termes $o(1)$ et $O(1)$ 
d\'ependant aussi du nombre $m$ de fonctions.

\subsection{Une r\'ecurrence}
Soit maintenant $f$ une fonction avec $0\leq f\leq\nu$. Nous allons
utiliser la construction pr\'ec\'edente de fa\c con r\'ep\'et\'ee, construisant de 
proche en proche une suite $(f_j)$ de fonctions, une suite $(\CB_j)$
de $\sigma$-alg\`ebres et une suite $(\Omega_j)$ d'ensembles appartenant
\`a $\CB_j$. 

Posons $f_1=f$ et soient $\CB_1$ la $\sigma$-alg\`ebre grossi\`ere
$\{\emptyset,\Z_N\}$ et $\Omega_1=\emptyset$. 

Supposons que les constructions ont \'et\'e faites jusqu'au rang $j$.
Nous posons
$$
f_{j+1}=(1-\one _{\Omega_j})\,\bigl(f-\E(f\mid\CB_j)\bigr)
$$
et distinguons deux cas. 
\begin{itemize}
\item
Si $\nnorm{f_{j+1}}k\leq\varepsilon$
nous arr\^etons l'algorithme; 
\item
sinon nous  employons la m\'ethode pr\'ec\'edente avec la famille de
fonctions $(f_1,\dots,f_{j+1})$ pour d\'efinir la
$\sigma$-alg\`ebre $\CB_{j+1}$ et l'ensemble $\Omega_{j+1}$ et nous
it\'erons l'algorithme, ce qui est possible car $|f_{j+1}|\leq
1+\nu$, \`a multiplication pr\`es par un terme de la forme
$1+\sigma^{1/2}\,O_\varepsilon(1)$.
\end{itemize}
Montrons que cet algorithme s'arr\^ete apr\`es un nombre born\'e d'\'etapes.
S'il ne s'arr\^ete pas \`a l'\'etape $j$ alors 
$\dis \nnorm{f_{j+1}}k^{2^k}= \E(f_{j+1}\,\CD
f_{j+1})\geq\varepsilon^{2^k}\ .$
Imaginons pour simplifier que le petit ensemble $\Omega_j$ est vide. 
Nous aurions alors $f_{j+1}= f-\E(f\mid\CB_j)$ et comme par construction 
$\CD f_{j+1}$ est uniform\'ement proche de $\E(\CD f_{j+1}\mid\CB_{j+1})$, 
cette in\'egalit\'e nous permettrait de minorer 
\begin{multline*}
 \E\Bigl( \bigl( f-\E(f\mid\CB_j)\bigr)\cdot \E(\CD f_{j+1}\mid
\CB_{j+1})\Bigr)\\
 \E\Bigl( \bigl(\E(f\mid\CB_{j+1})-\E(f\mid\CB_j)\bigr)\cdot \E(\CD f_{j+1}\mid
\CB_{j+1})\Bigr)
\end{multline*}
et donc aussi $\norm{\E(f\mid\CB_{j+1})-\E(f\mid\CB_j)}_{L^2}$ en
appliquant l'in\'egalit\'e de Cauchy-Schwarz. 
Les calculs pr\'ecis permettent en fait de borner inf\'erieurement la quantit\'e
positive
$$
 \norm{(1-\one_{\Omega_{j+1}})\cdot\E(f\mid\CB_{j+1})}^2_{L^2}
- \norm{(1-\one_{\Omega_{j}})\cdot\E(f\mid\CB_{j})}^2_{L^2}\ .
$$
Comme $\norm{(1-\one_{\Omega_{j}})\cdot\E(f\mid\CB_{j})}_{L^2}$ est
major\'e  par $\norm f_{L^2}$, cela prouve que l'algorithme s'arr\^ete en
temps born\'e.

La derni\`{e}re $\sigma$-alg\`ebre et le dernier 
ensemble $\Omega$ construits v\'erifient alors les propri\'et\'es de la
proposition~\ref{prop:decomp}.\hfill\qed

\section{Un poids pseudo-al\'eatoire}
\label{sec:pseudo}
Il nous reste \`a construire une fonction $f$ nulle en dehors des nombres
premiers et un poids
pseudo-al\'eatoire $\nu$  v\'erifiant les hypoth\`eses du
th\'eor\`eme~\ref{th:GTSz}.

Dans cette section
 nous notons $P$ l'ensemble des nombres 
premiers et la lettre $p$ d\'esigne toujours un
nombre premier. Rappelons la d\'efinition de deux fonctions classiques.
\begin{itemize}
\item
$\phi$ est la fonction indicatrice d'Euler :
pour tout entier $x>0$,  $\phi(x)$ est le nombre d'entiers compris
entre $1$ et $x$ et premiers avec $x$. 
\item
$\mu$ est la fonction de
M\"obius :
$$
 \mu(x)=\begin{cases}
1 & \text{si $x=0$ ;}\\
(-1)^\ell & \text{si $x$ est le produit de $\ell$
nombres premiers distincts ;}\\
0 &\text{sinon.}
\end{cases}
$$
\end{itemize}
\subsection{La fonction $f$ et le poids $\nu$}
Nous nous donnons une fois pour toutes
\begin{itemize}
\item
 une fonction 
$\chi\colon\R\to\R^+$, de
classe $\CC^\infty$, \`a support dans $[-1,1]$, avec $\chi(0)>0$ et
$\int_0^\infty(\chi'(x))^2\,dx=1$ ;
\item
  $w=w(N)$ une fonction enti\`ere tendant vers l'infini tr\`es
lentement, par exemple de l'ordre de $\log\log N$.
\end{itemize}
De plus nous notons
\begin{itemize}
\item
$W=W(N)$ le produit des nombres premiers $\leq w(N)$ ;
\item
 $R=N^\alpha$ o\`u $\alpha$ est une constante positive suffisamment
petite.
\end{itemize}

Soit $b$ un entier qui sera d\'efini plus bas et d\'efinissons la fonction
$f$ sur $[0,N[$ par 
\begin{equation}
\label{eq:def-f}
f(n)=\begin{cases}
\frac{\phi(W)}{W}\,\log(Wn+b) & 
\text{si $ Wn+b\geq R$ et $Wn+b$ est premier}\\
0 & \text{sinon.}
\end{cases}
\end{equation}
Rappelons l'estimation \'el\'ementaire de  Tch\'ebytchev pour le nombre
$\pi(x)$ de nombres premiers inf\'erieurs o\`u \'egaux \`a $x$:
$$
 c_1\frac x{\log x}\leq \pi(x)\leq c_2\frac x{\log x}
$$
o\`u $c_1$ et $c_2$ sont des constantes positives.  
 On en d\'eduit (sans utiliser le th\'eor\`eme de Dirichlet) qu'on peut
choisir $b\in [1,W[$, premier avec $W$ et tel que pour $N$ assez grand
on ait 
$\E(f(n)\mid n\in[0,N[)\geq\delta$, o\`u $\delta>0$ est une
constante.

Pour tout $n$ entier nous posons 
\begin{equation}
\label{eq:def-lambda}
 \lambda(n)=\sum_{d|n}\mu(d)\,\chi\bigl(\frac{\log d}{\log R}\bigr)
\end{equation}
et nous d\'efinissons
\begin{equation}
\label{eqq:def-nu}
 \nu(n)=\frac{\phi(W)}W\,\log R\cdot \lambda^2(Wn+b)\ .
\end{equation}

On v\'erifie facilement que pour tout $n\in[0,N[$ on a $0\leq f(n)\leq
c\nu(n)$ pour une certaine constante positive $c$. Ainsi 
la fonction $f$ v\'erifie (\`a une normalisation pr\`es) 
les propri\'et\'es annonc\'ees dans le th\'eor\`eme~\ref{th:existe-nu}. 

Il reste
\`a montrer que  $\nu$ est un poids pseudo-al\'eatoire, c'est \`a dire
que cette fonction v\'erifie la condition sur les formes lin\'eaires et la condition 
sur les corr\'elations. Nous ne donnons pas ici la preuve de ces
propri\'et\'es et nous contentons de montrer comment obtenir l'estimation
$$\E(\nu\mid\Z_N)=1+o(1)\ .$$
La m\'ethode qui suit n'est sans doute pas la plus simple possible mais elle
 contient les principaux ingr\'edients utilis\'es
dans la d\'emonstration compl\`ete et \'eclaire le r\^ole jou\'e par le
param\`etre $R$  et par la fonction $\chi$.
Intuitivement, le r\^ole du param\`etre  $w$ est d'\'eliminer les
perturbations
produites par les petits nombres premiers. Soit en effet $p\in P$.
Les nombres premiers se r\'epartissent
dans $p-1$ classes de congruence modulo $p$, ce qui cause une
irr\'egularit\'e d'ordre $1/p$, non n\'egligeable si $p$ est trop petit devant $N$. 

\subsection{Une r\'e\'ecriture}

$\E(\nu\mid\Z_N)$ est le produit par $C=W^{-1}\phi(W)\log R$ de

\begin{equation}
\label{eq:sum-d}
\sum_{d,d'}\mu(d)\mu(d')\,\chi\bigl(\frac{\log d}{\log R}\bigr)\,
\chi\bigl(\frac{\log d'}{\log R}\bigr)\,\E(\one_{d,d'}(x)\mid x\in\Z_N)
\end{equation}
o\`u $d,d'$ sont des entiers positifs et 
$$
 \one_{d,d'}(x)=\begin{cases} 
1 & \text{si $\ppcm(d,d')$ divise $Wx+b$ ;}\\
0 & \text{sinon.}
\end{cases}
$$
Nous \'evaluons $\E(\one_{d,d'}\mid\Z_N)$. 
\`A cause des facteurs $\mu(d)$ et $\mu(d')$ nous pouvons nous
restreindre au cas o\`u $d$ et $d'$ sont sans carr\'e et \'ecrire
$$
 d=\prod_p p^{\omega_p}\ ;\  d'=\prod_p p^{\omega'_p}
$$
o\`u $\bomega=( \omega_p\;;\;p\in P)$ et  $ \bomega'=( \omega'_p\;;\;p\in P)$ sont des
suites \`a valeurs dans $\{0,1\}$.
D\'efinissons $E_p(\eta)$ pour $p\in P$ et $\eta\in\{0,1\}$ par
\begin{equation}
\label{eq:def-Ep}
E_p(0)=1\ ;\ E_p(1)=0\text{ si } p\leq w\text{ et }
E_p(1)=\frac 1p\text{ si }p>w\ .
\end{equation}
Nous avons
\begin{equation}
\label{eq:esp-one}
 \E(\one_{d,d'}\mid\Z_N)=\prod_p
E_p(\max( \omega_p, \omega'_p))+N^{-1}\,O(1)\ .
\end{equation}
En effet, 
s'il existe un premier $p$ avec $p\leq w$ et
$\max( \omega_p, \omega'_p)=1$ alors $\one_{d,d'}(x)=0$ pour tout $x$
puisque $p$ divise $W$ et est premier avec $b$. 
Dans le cas contraire la proportion des $x\in\Z_N$ tels que $\ppcm(d,d')$ divise $Wx+b$ 
est $1/\ppcm(d,d')$ \`a un $N^{-1}\,O(1)$ pr\`es. 

Nous comprenons maintenant le r\^ole de la troncature effectu\'ee par la fonction $\chi$ dans la
d\'efinition de $\nu$ : pour tous les $d$ consid\'er\'es, les entiers
$1,\dots,N-1$ se r\'epartissent suffisamment uniform\'ement dans les
classes de congruence modulo $d$. En effet, la 
somme~\eqref{eq:sum-d} contient au maximum $R^2$ termes non nuls  et
la somme des erreurs  $N^{-1} \,O(1)$  est donc de la forme
$R^2N^{-1}O(1)$ ;
le produit par $C$ de cette 
expression est un $o(1)$ puisque $R$ est une petite puissance de $N$.

Ainsi, en reportant l'estimation~\eqref{eq:esp-one} dans~\eqref{eq:sum-d} et 
en rempla\c cant la fonction de M\"obius par sa d\'efinition, nous
obtenons
\begin{multline}
\label{eq:approx}
\E(\nu\mid\Z_N)
=o(1)\\ +C
\sum_{ \bomega, \bomega'}
\chi\Bigl(\frac{\sum_p { \omega_p}\log p}{ \log R}\Bigr)\,
\chi\Bigl(\frac{\sum_p { \omega'_p}\log p}{ \log R}\Bigr)\,
\prod_p(-1)^{ \omega_p+ \omega'_p}
E_p(\max( \omega_p, \omega'_p))\; .
\end{multline} 

\subsection{Transform\'ee de Fourier}
Dans l'article~\cite{GrT}, Green et Tao utilisent une troncature brutale de
la somme~\eqref{eq:def-lambda} d\'efinissant la fonction $\lambda$, comme dans les travaux
de Goldston et Y\i ld\i r\i m. L'emploi de la fonction lisse $\chi$ 
permet de
se servir de la transform\'ee de Fourier. Nous \'ecrivons
\begin{equation}
\label{eq:fourier}
 \chi(x)=\int \tau(t) e^{-x(1+it)}\,dt
\end{equation}
o\`u $\tau$ est une fonction \`a d\'ecroissance rapide. La
somme dans la formule~\eqref{eq:approx} se met alors ais\'ement sous la forme
\begin{equation}
\label{eq:prod-sum}
\iint \tau(t)\tau(t')\,\prod_{p} Z_p(t,t')\,dt\,dt'
\end{equation}
o\`u, en notant 
$$
 z=\frac{1+it}{\log R} \text{ et } z'=\frac{1+it'}{\log
R}
$$
nous avons
$$
 Z_p(t,t')=\sum_{\eta,\eta'\in\{0,1\}} (-1)^{\eta+\eta'}
\,E_p(\max(\eta,\eta'))\,
p^{-\eta z-\eta' z'}\ .
$$
En rempla\c cant les $E_p$ par leurs valeurs~\eqref{eq:def-Ep} nous
obtenons
$$ Z_p(t,t')=\begin{cases}
1 &\text{si $p\leq w$ ;}\\
1-p^{-1-z}-p^{-1-z'}+p^{-1-z-z'} & \text{ si $p>w$.}
\end{cases}
$$
Il vient alors
\begin{multline*}
 \prod_p
Z_p(t,t')=\prod_{p>w}\bigl(1-p^{-1-z}-p^{-1-z'}+p^{-1-z-z'}\bigr)\\
=(1+o(1))\prod_{p>w} \frac{(1-p^{-1-z})(1-p^{-1-z'})}{1-p^{-1-z-z'}}
\end{multline*}
car $w$ tend vers l'infini avec $N$
 et donc
\begin{equation}
\label{eq:prodZp}
\prod_{p}Z_p(t,t')=(1+o(1))\frac{\zeta(1+z+z')}{\zeta(1+z)\zeta(1+z')}\,
\Bigl(\prod_{p\leq w}
\frac{(1-p^{-1-z})(1-p^{-1-z'})}{1-p^{-1-z-z'}}\Bigr)^{-1}\ .
\end{equation}
Nous \'ecrivons l'int\'egrale~\eqref{eq:prod-sum} comme somme de
l'int\'egrale pour $t,t'$ appartenant \`a l'intervalle
 $J=[-(\log R)^{1/2},(\log R)^{1/2}]$ et d'un reste qui est de la
forme $C^{-1}o(1)$ \`a cause de la d\'ecroissance rapide de la fonction
$\tau$.

Soient maintenant $t$ et $t'$ appartenant \`a $J$.\\
Comme $W$ est tr\`es petit
devant $(\log R)^{1/2}$ nous avons
$$
\prod_{p\leq w}\frac{(1-p^{-1-z})(1-p^{-1-z'})}{1-p^{-1-z-z'}} =
(1+o(1))\, \prod_{p\leq w}(1-p^{-1})=(1+o(1))\, \frac{\phi(W)}W\ .
$$
D'autre part, l'estimation \'el\'ementaire
$$
 \zeta(1+s)\sim \frac 1s\text{ quand }s\to 0\text{ avec }\reel(s)>0
$$
nous donne
$$
\frac{\zeta(1+z+z')}{\zeta(1+z)\zeta(1+z')}=(1+o(1))\,
\frac{zz'}{z+z'}= (1+o(1))\,\frac 1{\log R}\,\frac{(1+it)(1+it')}{2+it+it'}\,
\ .
$$
En reportant ces valeurs dans~\eqref{eq:prodZp} nous obtenons
$$
 \iint_{J\times J}\prod_p Z_p(t,t')=C^{-1}\,(1+o(1))\,\iint_{J\times J}
\tau(t)\tau(t')\,\frac{(1+it)(1+it')}{2+it+it'}\,dt\,dt'\ .
$$
Comme $\tau$ est \`a d\'ecroissance rapide, l'int\'egrale de la
deuxi\`eme fonction  
 en dehors de $J\times J$ est de un $C^{-1}\,o(1)$ et nous obtenons
$$
 \iint\prod_p Z_p(t,t')
=C^{-1}\,(1+o(1))\,\iint
\tau(t)\tau(t')\frac{(1+it)(1+it')}{2+it+it'}\,dt\,dt'\ .
$$
Cette derni\`ere int\'egrale s'\'ecrit
$$
 \int_0^{+\infty}\,ds\iint
\tau(t)\tau(t')(1+it)(1+it')e^{-s(2+it+it')}\,dt\,dt'
=\int_0^{+\infty}(\chi'(s))^2\,ds= 1
$$
et on a donc bien $\E(\nu\mid\Z_N)=1+o(1)$.\hfill\qed

\end{document}